\documentclass[12pt]{amsart}
\usepackage{amsmath,amssymb,amscd,amsthm}
\usepackage{euscript}
\begin{document}
\renewcommand{\theequation}{\thesection.\arabic{equation}}

\newcommand{\supp}{\operatorname{supp}}
\newcommand{\dist}{\operatorname{dist}}

\newtheorem{thm}{Theorem}
\newtheorem{cas}{Case}

\theoremstyle{definition}
\newtheorem*{rem}{Remark}
\newtheorem*{rems}{Remarks}
\newtheorem*{claim}{Claim}

\title{The Calder\'on Commutator along a
Parabola}
\thanks{   
The authors were
supported in part by EPSRC. Hofmann
supported in part by an NSF grant. Research
at MSRI supported in part by NSF grant DMS-9701755.}

\author{Anthony Carbery}
\address{\hskip-\parindent Anthony Carbery\\
Department
of Mathematics and Statistics\\ University of
Edinburgh\\ Edinburgh\\Scotland EH9 3JZ}
\email{carbery@maths.ed.ac.uk}

\author{Steve Hofmann}
\address{\hskip-\parindent Steve Hofmann\\Department of
Mathematics\\ University of Missouri\\
Columbia, MO 65211\\ USA}
\email{hofmann@math.missouri.edu}

\author{James Wright}
\address{\hskip-\parindent James Wright\\
Department of Mathematics\\
University of New South Wales\\ Sydney,
N.S.W. 2052\\ Australia}
\email{jimw@maths.unsw.edu.au}

\begin{abstract} We introduce an
analogue of Calder\'on's first
commutator along a parabola, and
establish its
$L^2$ boundedness under essentially
sharp hypotheses.\end{abstract}
\maketitle

\section{Introduction}

During the past 30 years, many authors
have investigated the $L^p$ mapping
properties of singular integral
operators whose ``kernels'' are actually
singular measures supported on lower
dimensional subvarieties. One may
consult the survey article of Wainger
\cite{W} for a history of the subject up
to the mid 1980's. The interested reader
may also find a synopsis of more recent
developments, along with numerous
references, in the monograph of Stein
\cite[Chapter XI, sections 4.5, 4.7,
4.17]{St}. The prototypical example of
of such operators is the ``Hilbert
transform along a curve'':
\begin{equation} H_\gamma f(x)\equiv
p.v.\int_{\mathbb{R}}f(x-\gamma
(t))\frac{dt}{t},
\label{eq1.1}\end{equation} where $x\in
\mathbb{R}^n$ and $\gamma :\mathbb{R}\to
\mathbb{R}^n$ is the parametrization of
a smooth curve in
$\mathbb{R}^n$. This sort of operator
was first introduced by Fabes \cite{F}
who established
$L^2$ boundedness of
$H_\gamma$ in the special case $n=2$,
$\gamma(t)=(t,(sgnt)t^2)$. The original
motivation for studying operators of
this sort is that they arise when one
tries to develop a ``method of
rotations'' (see \cite{CZ} in the
classical elliptic case) for parabolic
singular integrals of convolution type.
In turn, the method of rotations enables
one to significantly relax the
regularity hypotheses on kernels of
singular integral operators.

In the present paper, motivated in part
by formal analogy to Fabes' goal of
extending the method of rotations to the
parabolic case, and in part by recent
developments in the theory  and
application of parabolic singular
integrals which are not of convolution
type \cite{LM},
\cite{H1}, \cite{H2}, \cite{HL}, we
introduce a certain bilinear analogue of
\eqref{eq1.1}: the ``Calder\'on
commutator along a parabola''. To
describe this operator, suppose that
$A:\mathbb{R}^2\to \mathbb{R}$ satisfies
the $Lip_{1,\frac{1}{2}}$ condition
\begin{equation}\label{eq1.2}
|A(x+h)-A(x)|\leq B_0r,\end{equation}
for some constant $B_0$, whenever
$h\equiv (h_1,h_2)$ satisfies $|h_1|\leq
r$, $|h_2|\leq r^2$. We remark that the
results of the present paper may be
readily extended to $\mathbb{R}^n$,
$n>2$, by the same arguments which we
shall give below, but to minimize
technicalities we shall restrict our
attention to the case
$n=2$. Since the story here is not yet
complete (we do not yet know how to
treat the higher order commutators, for
example), it does not seem crucial at
this point that we state our results in
the greatest generality.

For $A$ as in  \eqref{eq1.2} and for
$\gamma(t)\equiv (t,t^2)$, we define
\begin{equation} T^\gamma _Af(x)\equiv
p.v.\int_\mathbb{R}[A(x)-A(x-
\gamma(t))]f(x-\gamma
(t))\frac{dt}{t^2}.
\label{eq1.3}\end{equation} As mentioned
above, the operator $T_A$ has a
connection with the results of
\cite{LM}, \cite{H1}, \cite{H2}, and
\cite{HL}, which may be understood as
follows. Let $K:\mathbb{R}^2\backslash
\{0\}\to
\mathbb{R}$ satisfy the parabolic
homogeneity property
\begin{equation}\label{eq1.4}K(\rho
x_1,\rho ^2x_2)\equiv
\rho^{-4}K(x),\end{equation} and further
suppose that
$K\in L^1(S^1)$, and that $K(x_1,x_2)$
is even in $x_1$, for each fixed $x_2$.
Then one may use parabolic polar
co-ordinates to obtain a representation
of the parabolic Calder\'on commutator
in terms of $T^\gamma_A$. Indeed,
\begin{equation}\begin{split}
\label{eq1.5} C_Af(x)&\equiv
p.v.\int_{\mathbb{R}^2}[A(x)-
A(y)]K(x-y)f(y)dy\\ &\equiv p.v.
\int_{S^1}K(\sigma
)T^{\gamma_\sigma}_Af(x)(1+
\sigma^2_2)d\sigma,
\end{split}\end{equation} where
$T^{\gamma_\sigma}_A$ is defined as in
\eqref{eq1.3}, with $\gamma_\sigma
(t)\equiv (t\sigma _1,t^2\sigma_2)$,
$(\sigma_1,\sigma_2)\in S^1$. We have
used here the parabolic polar coordinates
\begin{align*}x_1&=\rho \sigma_1,\quad
x_2=\rho ^2\sigma_2\\
dx&=\rho^2d\rho(1+\sigma_2^2)d\sigma.
\end{align*} When $\sigma=(\pm 1,0)$,
$L^2$ boundedness of
$T^{\gamma_\sigma}_A$ reduces to that of
Calder\'on's original commutator on the
line
\cite{C}; when $\sigma =(0,\pm 1)$,
matters reduce to a result of Murray
\cite{Mu}. Otherwise,
$L^2$ boundedness of
$T^{\gamma_\sigma}_{A}$ is new and will
be treated in this paper, although to
simplify the notation we shall take
$\sigma_1=1=\sigma_2$. The proof in the
general case is identical and one
obtains bounds independent of
$\sigma \in S^1$.

In \cite{H1}, it was
shown that for the special case
\begin{equation}\label{eq1.6}
K(x_1,x_2)\equiv x_2^{-2}\exp
\left\{\frac{-|x_1|^2}{4x_2}\right\}
\chi _{\{x_2>0\}},\end{equation} a
necessary and sufficient condition for
the $L^2$ boundedness of $C_A$ is that,
for some $B<\infty$,
\begin{equation}\begin{split}\text{(i)}
\quad &|A(x_1+h,x_2)-A(x_1,x_2)| \leq
B|h|\\
\text{(ii)}\quad
&\|\mathbb{D}_2A\|_{BMO}\leq
B,\label{eq1.7}\end{split}\end{equation}
where
$$\mathbb{D}_2A\equiv
\left(\frac{\xi_2}{\sqrt{|\xi
_1|^2-i\xi_2}}\hat{A}(\xi
)\right)^\vee,$$ and $\hat{f}$ and $f^\vee$ denote
the Fourier and inverse Fourier transforms of $f,$
respectively. The $BMO$ norm is defined as usual
by
$$ \|b\|_{BMO} = sup \int_I |b - m_I b| ,$$
where $m_I b$ denotes the mean value of $b$ over $I,$ and
where the sup runs over all parabolic
``cubes'' of the form
\begin{equation}\label{eq1.9}
I \equiv I_r(a,b)\equiv
\left[a-\frac{r}{2},a+\frac{r}{2}
\right]\times
\left[b-\frac{r^2}{2},b+
\frac{r^2}{2}\right].\end{equation} The
higher dimensional case was also treated
in
\cite{H1}. The kernel \eqref{eq1.6}
arises in the multilinear expansion of
the parabolic double layer potential on
a time varying graph $x_3=A(x_1,x_2)$. A
condition similar to \eqref{eq1.7}
(later shown to be equivalent in
\cite{HL}) had previously been shown to
be sufficient for
$L^2$ boundedness in \cite{LM}. It was
also shown in \cite{H1} that
\eqref{eq1.7} implies the
$Lip_{1,1/2}$ condition \eqref{eq1.2},
with bound $B_0\equiv CB$.

In this paper we prove the following:

\setcounter{thm}{0}\begin{thm}\label{t1}
Suppose $A$ satisfies \eqref{eq1.7}. Then
$T^\gamma_A$ is bounded on
$L^2$.\end{thm}

\begin{rems}
\hspace*{.2in}
\begin{itemize}
\item[(1)]
By our previous comments, one
may then obtain a ``method of
rotations'' which permits us to deduce
$L^2$ boundedness of $C_A$ (defined in
\eqref{eq1.5}) for $K$ as in
\eqref{eq1.4} satisfying only that $K\in
L^1(S^1)$, and that $K(x_1,x_2)$ is even
in $x_1$, for each fixed $x_2$.
Previously, it had been required that
$K$ have some smoothness on
$S^1$.
\item[(2)]
One can also show that $T^\gamma_A$ is bounded on 
$L^p , 1<p<\infty$, by using the key estimate
\eqref{eq2.4} below and Littlewood-Paley arguments
adapted to rough singular integral operators. See e.g.,
\cite{H3}.
\end{itemize}
\end{rems}

In the next section we begin the proof
of Theorem~\ref{t1}, and complete it in Section 3.

\section{Proof of Theorem~\ref{t1}}
\setcounter{equation}{0}

We may assume, and do, that the constant
$B$ in \eqref{eq1.7} is one. By the
parabolic version of \cite{H3} (whose
proof is virtually identical that given
in \cite{H3}, and is therefore omitted
here), it is enough to check the
following ``rough operator'' $T1$
criterion, which consists of three
parts. We need to prove that
\begin{equation}\label{eq2.1}
T^\gamma_A1,
\left(T^\gamma_A\right)^*1\in
BMO,\end{equation} where
$\left(T^\gamma_A\right)^*$ is defined
as a mapping from test functions to
distributions by
$$\langle
\left(T^\gamma_A\right)^*f,g\rangle
\equiv
\langle f,T^\gamma_Ag\rangle,$$ for all
$f,g\in C^\infty_0(\mathbb{R}^2)$.
Second, for all $x\in \mathbb{R}^2$,
$r>0$, let $\Phi (x,r)$ denote the class
of all
$\varphi \in C^\infty_0$, supported in
$I_r(x)$ (see \eqref{eq1.7}), and
satisfying
\begin{align}\label{eq2.2}&\text{(i)}&
\|\varphi \|_\infty &\leq 1\\
&\text{(ii)}& |\varphi (y+h)-\varphi
(y)|&\leq \rho /r,\notag\\
\intertext{whenever $|h_1|\leq\rho$,
$|h_2|\leq \rho^2$, $h\equiv (h_1,h_2)$,
and} &\text{(iii)}&\sup_{k+m\leq
2}\|\left(\frac{\partial}{\partial
x_1}\right)^k\left(
\frac{\partial}{\partial
x_2}\right)^m\varphi \|_\infty &\leq
1.\notag\end{align} We shall need to
establish the Weak Boundedness Property
(WBP):
\begin{equation}\label{eq2.3} |\langle
\psi ,T^\gamma_A\phi \rangle |\leq
Cr^3,\end{equation} for all $\psi,\phi
\in \Phi (x,r)$, for any $r>0$, $x\in
\mathbb{R}^2$. Since the homogeneous
dimension of parabolic $\mathbb{R}^2$ is
$d=3$, \eqref{eq2.1} and \eqref{eq2.3}
are the usual David-Journ\'e conditions
\cite{DJ} in this context. However, in
lieu of the standard Calder\'on-Zygmund
kernel conditons considered in
\cite{DJ},we shall establish instead the
weak smoothness condition of \cite{H3},
which had also appeared in a similar connection in
\cite{CS}.
Let $Q_s$ denote the operator defined by
convolution with a smooth function
$\psi_s$ which has mean value zero, is
supported in $I_s(0)$, and is normalized
so that $\|\psi_s\|_1=1$. We shall prove
that, whenever $s\leq 2^j$,
\begin{equation}\label{eq2.4}\|Q_sT_j
\|_{L^2\to L^2}\leq
C\left(\frac{s}{2^j}\right)^\epsilon
\end{equation} for some $\epsilon >0$,
where
$$T_jf(x)\equiv
\int_\mathbb{R}[A(x)-A(x-\gamma
(t))]f(x-\gamma (t))\eta
\left(\frac{t}{2^j}\right)
\frac{dt}{t^2},$$ and where $\eta \in
C^\infty_0\left[\left(\frac{1}{2},2
\right)\cup
\left(-2,-\frac{1}{2}\right)\right]$,
$0\leq
\eta \leq1$, and
$\sum^\infty_{j=-\infty}\eta
\left(\frac{\cdot}{2^j}\right)\equiv 1$
away from $0$.

By \cite{H3} (or rather its parabolic
analogue), Theorem~\ref{t1} follows
immediately from
\eqref{eq2.1},
\eqref{eq2.3} and \eqref{eq2.4}. In this
section we shall prove \eqref{eq2.1} and
\eqref{eq2.3}. In the next section we
prove \eqref{eq2.4}.

To establish WBP \eqref{eq2.3}, by
dilation invariance we may assume $r=1$,
so it is clearly enough to prove that,
for any $x_0\in \mathbb{R}^2$,
\begin{equation}\label{eq2.5}
T^\gamma_A\varphi \in L^2(I)+L^\infty
(I),\quad \forall \varphi
\in \Phi (x_0,4),\end{equation} where
$I=I_1(x_0)$. We claim that
\eqref{eq2.5} also implies that
$T^\gamma_A1\in BMO$ (we omit the proof
of the fact that $(T^\gamma_A)^*1\in
BMO$, since it is identical). To prove
the claim, we need to consider
\begin{equation}\label{eq2.6}
\frac{1}{|I|}\int_I|T^\gamma_A
1-C_I|,\end{equation} where
$I=I_r(x_0)$, and by dilation invariance
we may take $r=1$. We write
$$1=\phi +(1-\phi ),$$ where $\phi
\equiv 1$ in $I_3(x_0)$, and $\phi \in
C^\infty_0(I_4(x_0))$. Then up to a
normalizing constant,
$\phi \in \Phi (x_2,4)$, so by
\eqref{eq2.5}
$$\int_I|T^\gamma_A\phi |\leq C.$$ We
now set $C_I\equiv T^\gamma_A(1-\varphi
)(x_0)$, so that, for all $x\in I$,
\begin{equation*}\begin{split}&|
T_A^\gamma (1-\phi
)(x)-T^\gamma_A(1-\phi )(x_0)|\\ &\quad
\leq\int|A(x)-A(x-\gamma
(t))-[A(x_0)-A(x_0)-\gamma (t)]||1-\phi
(x-\gamma (t))|\frac{dt}{t^2}\\ &\quad
+\int|A(x_0)-A(x_0-\gamma(t))||\phi
(x-\gamma (t)-\phi (x_0-\gamma
(t))|\frac{dt}{t^2}\\ &\quad \leq
c\left(\int_{|t|>1}\frac{dt}{t^2}+
\int_{1<|t|<5}\frac{dt}{t}\right)=
C,\end{split}\end{equation*} by
\eqref{eq1.2} (with $B_0\equiv C$), and
the definition of $\phi $. Consequently,
\eqref{eq2.6} is no larger than $C$,
which proves the claim. Thus, we have
reduced the proofs of \eqref{eq2.1} and
\eqref{eq2.3} to that of
\eqref{eq2.5}.

We establish \eqref{eq2.5} under the a
priori assumption that $A\in C^\infty$,
but our quantitative estimates will
depend only on the bounds in
\eqref{eq1.7} (which we have normalized
to be 1). Let $\phi
\in
\Phi (x_0,4)$ and then integrate by
parts; to do this rigorously requires
truncation of the principal value
integral in \eqref{eq1.3}, but it is
routine to verify that the boundary
terms (which arise when integrating the
truncated integrals by parts) are
harmless. We shall therefore argue
formally, and ignore all such
truncations and boundary terms, in order
not to tire the reader with minutiae.
Formally then,
\begin{equation*}\begin{split}
T^\gamma_A\phi
(x)&=\int_\mathbb{R}\frac{\partial
A}{\partial x}(x-\gamma (t))\phi
(x-\gamma (t))\frac{dt}{t}
+2\int_\mathbb{R}\frac{\partial
A}{\partial t}(x-\gamma (t))\phi
(x-\gamma (t))dt\\ &+\int[A(x)-
A(x-\gamma (t)]\left\{\frac{\partial
\phi}{\partial x}(x-\gamma
(t))+2\frac{\partial \phi }{\partial
t}(x-\gamma
(t))t\right\}\frac{dt}{t}\equiv
I+II+III.\end{split}\end{equation*} But
$I$ is precisely the Hilbert transform
along the curve $(t,t^2)$, acting on the
$L^2$ function
$\frac{\partial A}{\partial x}\phi $.
Thus, $I\in L^2$. Also $III\in L^\infty
(I)$, by
\eqref{eq1.2} and the fact that $\phi
\in \Phi (x_0,4)$ (see 2.2). The only
delicate term is II. We recall that
$\frac{\partial A}{\partial x}\in
L^\infty, \mathbb{D}_2A\in BMO$, with
norm $B=1$ (given our normalization).
Define $a\equiv
\mathbb{D}A\equiv
(\sqrt{|\xi_1|^2-i\xi_2}\hat{A}(\xi
))^\vee$, which is therefore in $BMO$,
with norm $CB=C$ by parabolic
Calder\'on-Zygmund Theory
\cite{FR1,FR2}. Then
$$\frac{\partial A}{\partial t}\equiv
\mathbb{D}_2(\mathbb{D}A)\equiv
\mathbb{D}_2a,$$ so that
\begin{equation}\label{eq2.7}II=
\int\mathbb{D}_2a(x-\gamma (t))\phi
(x-\gamma (t))dt.\end{equation} By
\cite[Lemma 1 and Lemma 2, pp.
111-113]{FR1},
$$\mathbb{D}_2a(x)\equiv p.v.\int
k(x-y)a(y)dy,$$ where $k(x)$ is odd,
belongs to $C^\infty
(\mathbb{R}^2\backslash \{0\})$ and
satisfies the homogeneity property
\begin{equation}\label{eq2.8}k(\lambda
x,\lambda^2t)\equiv
\lambda^{-d-1}k(x,t)\end{equation} (we
recall that
$d=3$ is the homogeneous dimension of
parabolic $\mathbb{R}^2$). By the
oddness of $k$, we may assume that a has
mean value zero on $I_1(x_0)\equiv I$.
Let $\|x\|\equiv |x_1|+|x_2|^{1/2}$ be
the parabolic ``norm'' of $x$. Let
$\eta \in C^\infty_0[-10,10]$, $\eta
\equiv 1$ on $[-9,9]$, and set
$a_1(x)=a(x)\eta (\|x-x_0\|)$,
$a_2=a-a_1$. By \eqref{eq2.8} and the
parabolic version of a standard estimate
of \cite{FS},
\begin{equation*}\begin{split}|
\mathbb{D}_2a_2(x)|&\leq
\int\frac{c}{1+\|x-y\|^{d+1}}|a(y)|dy\\
&\leq C\|a\|_{BMO}
\leq C,\end{split}\end{equation*}
whenever $x\in \supp \phi \subseteq
I_4(x_0)$. Consequently the contribution
of
$a_2$ to \eqref{eq2.7} yields a bounded
term, as desired. The contribution of
$a_1$ is
\begin{equation*}\begin{split}&\int\phi
(x-\gamma (t))\eta
(t)\mathbb{D}_2a_1(x-\gamma (t))dt\\
\intertext{(since $x\in I_1(x_0)$,
$x-\gamma (t)\in I_4(x_0)$)} &\quad
=\int\eta (t)\int k(x-\gamma (t)-y)[\phi
(x-\gamma (t))-\phi (y )]a_1(y)dydt\\
&\quad\quad + \int \eta
(t)\mathbb{D}_2(\phi a_1)(x-\gamma
(t))dt\\ &\quad
\equiv
II_1+II_2.\end{split}\end{equation*} By
\cite{H1} (or even \cite{FR1}, since
$\phi $ is smooth), the commutator
$[\mathbb{D}_2,\phi ]$ defines a bounded
operator on $L^2$. Hence, by Minkowksi's
inequality,
\begin{equation*}\begin{split}\|II_1
\|_2&\leq
\int\eta (t)\|[\mathbb{D}_2,\phi
]a_1(\cdot -\gamma (t) \|_2dt\\ &\leq
c\|a_1\|_2\leq c\| a\|_{BMO}\leq
C.\end{split}\end{equation*} Finally,
since $\phi a_1\in L^2$, and
$\mathbb{D}_2$ is given by the multiplier
$\xi_2(|\xi_1|^2-i\xi_2)^{-1/2}$, it is
enough to show that the operator
$$f\to \int f(x-\gamma (t))\eta (t)dt$$
is smoothing on $L^2$ of parabolic order
$1$; i.e., that the multiplier
\begin{equation*}\begin{split}m(\xi
)&\equiv \int e^{i(\xi_1t+\xi_2t^2)}\eta
(t)dt\\
\intertext{satisfies} |m(\xi )|&\leq
C(|\xi
_1|+|\xi_2|^{1/2})^{-1}
\end{split}\end{equation*} This estimate
is well known, and is
an easy consequence of standard
integration by parts arguments as may be found in
\cite[Chapter VIII]{St}. For the sake of
completeness, we sketch the argument here.

\begin{cas}\label{cas1} $|\xi _1|>40|\xi
_2|$.\end{cas} In this case,
\begin{equation*} |n(\xi )|=\left|\int
e^{i(\xi_1t+\xi_2t^2)}\frac{d}{dt}
\left(\frac{1}{\xi_1+2+\xi_2}\eta
^0(t)\right)dt\right|\leq
C|\xi_1|^{-1},\end{equation*} Since
$\supp \eta \leq [-10,10]$.

\begin{cas}\label{cas2}
$|\xi_1|<40|\xi_2|$.\end{cas} In this
case we seek the estimate
$|m(\xi )|\leq C|\xi _2|^{-1/2}$. We
write
$$m(\xi )=\int e^{i\lambda \varphi_\xi
(t)}\eta (t)dt,$$ where $\lambda \equiv
|\xi_2|$, and
$\phi _\xi (t)\equiv
\left(\frac{\xi_1}{|\xi _2|}t+sgn
\xi_2t^2\right)$. Note that $|\phi ''
_\xi (t)|=2$. The desired bound now
follows immediately from Vander Corput's
Lemma, or to be more precise, its
corollary given in \cite[p. 334,
inequality (6)]{St}. This concludes the
proof of
\eqref{eq2.5}, and therefore also the
proofs of \eqref{eq2.1} and
\eqref{eq2.3}. We finish the proof of
our Theorem in the next section, in
which we prove \eqref{eq2.4}.

\section{Proof of Theorem 1 (continued): estimate \eqref{eq2.4}}
\setcounter{equation}{0}

In this section, we give the proof of estimate (2.4), which
will complete the proof of Theorem 1.  This is the
most technical part of the proof, and follows ideas from \cite{Ch};
see also \cite{CVWWW} and \cite{CWW} .

By dilation invariance, we may take
$j=0$. We recall that $T_0f$ is a sum of
two terms
$$\tilde{T}_0f(x)\equiv
\int_\mathbb{R}[A(x)-A(x-\gamma
(t))]f(x-\gamma (t))\eta
(t)\frac{dt}{t^2}$$ where $\eta
$ is an even smooth cut-off function
supported in the ``half annulus''
$\frac{1}{2}<t<2$, plus  another term
with
$-2<t<-\frac{1}{2}$. By symmetry, we
treat only the former. Let $K_0$ denote
the kernel of
$\tilde{T}_0^* \tilde{T}_0^*$, where
$\tilde{T}_0^*$ denotes the adjoint of
$\tilde{T}_*$. Following \cite{Ch}, we see that the
desired bound \eqref{eq2.4} then follows
easily (we omit the routine details) from

\begin{claim}
\begin{equation}\label{eq3.1}
\int|K_0(x+h,y)-K(x,y)|dy\leq
C\lambda^{1/3},\end{equation} whenever
$|h_1|\leq
\lambda$, $|h_2|\leq \lambda^2$,
$\lambda\leq 1/1000$. Now, it is a
routine matter to see that
$$\tilde{T}_0^*g(x)=-\int[A(x)-A(x+
\gamma (s))]g(x+\gamma (s))\eta
(s)\frac{ds}{s^2}$$ Thus,

\begin{equation}\begin{split}
\label{eq3.2}
\tilde{T}_0\tilde{T}_0^*f(w)&=
\int[A(w)-A(w-\gamma (t))]T^*_0f(\omega
-\gamma (t))\phi(t)dt\\ &=-
\iint[A(\omega )-A(\omega -\gamma
(t))]\\ &\quad [A(\omega -\gamma
(t))-A(\omega +\gamma (s)-\gamma
(t)]f(\omega +\gamma (s)-\gamma (t))\phi
(t)\phi
(s)dtds,\end{split}\end{equation} where
for convenience of notation we have
defined

\begin{equation*}\phi (t)\equiv \eta
(t)t^{-2}.\end{equation*} Now, the claim
\eqref{eq3.1} amounts to saying that the
operator

\begin{equation*}f\to \int[K_0(\cdot
+h,y)-K_0( \cdot
,y)]f(y)dy\end{equation*} maps $L^\infty
\to L^\infty$ with norm
$C\lambda^{1/3}$. Hence, we may with
impunity excise (in
\eqref{eq3.2}) any set in $(s,t)$ space
of measure no larger than
$C\lambda^{1/3}$. Indeed, we shall
restrict the domain of integration in
\eqref{eq3.2} to the set
$$\left\{(s,t)\in
\left[\frac{1}{2},2\right]^2:|s-t|\geq
15\lambda ^{\frac{1}{3}}\right\}.$$ By
symmetry it suffices to consider only
the case $s>t$. Note that in the set
where $s-t\geq
\lambda ^{\frac{1}{3}}$, the map

\begin{equation}\begin{split}
\label{eq3.3}
\Phi _\omega :(s,t)&\mapsto \omega
+\gamma (s)-\gamma (t)=y\\
\frac{1}{4}\leq s<t\leq 4&\mapsto
\mathbb{R}^2\end{split}\end{equation} is
injective, and that the Jacobian matrix
is
$$\left[\begin{array}{cc} 1 & -1\\ 2s &
-2t\end{array}\right]$$ with determinant
$2(s-t)\geq 2\lambda ^{1/3}$. Let $S_0$
denote the part of the operator
$\tilde{T}_0\tilde{T}_0^*$ in
\eqref{eq3.2}, with $(s-t)>15\lambda
^{1/3}$. We then have that

\begin{equation*}\begin{split}|S_0
f(x+h)-S_0f(x)|&=|\iint_{E\lambda}\phi
(s)\phi (t)\left\{[A(x+h)-A(x+h-\gamma
(t)][A(x+h-\gamma (t)]\right.\\ &\quad
[A(x+h-\gamma (t))-A(x+h+\gamma
(s)-\gamma (t))] f(x+h+\gamma (s)-\gamma
(t))\\ &-[A(x)-A(x-\gamma (t))]\\ &\quad
\left.[A(x-\gamma (t))-A(x+\gamma
(s)-\gamma (t)]f(x+\gamma (s)-\gamma
(t))\right\}dsdt|\\ &\equiv
|\iint_{E_\lambda} \phi (s) \phi
(t)\left\{B(x+h,s,t)f(x+h+\gamma
(s)-\gamma (t))\right.\\
&-\left.B(x,s,t)f(x+\gamma (s)-\gamma
(t)\right\}dsdt|
\end{split}\end{equation*} (where
$E_\lambda\equiv \{(s,t): s-t\geq
15\lambda ^{1/3},\frac{1}{2}\leq t<s\leq
2\}$, and
$B(x,s,t)\equiv [A(x)-A(x-\gamma
(t))][A(x-\gamma (t))-A(x+\gamma
(s)-\gamma (t))]$)

\begin{equation*}\begin{split}
&=|\iint_{E_\lambda}\phi (s)\phi
(t)\left\{[B(x+h,s,t)-B(x,s,t)]f
(x+h+\gamma (s)-
\gamma (t))\right.\\ &\quad
+B(x,s,t)[f(x+h+\gamma (s)-\gamma
(t))-f(x+\gamma (s)-\gamma (t))]dsdt|\\
&\equiv |\iint
(I+II)dsdt|.\end{split}\end{equation*}
Now suppose that
$\|f\|_{L^\infty=1}$. Taking the
supremum over all such $f$, we see that
the contribution of
$\iint Idsdt$ to \eqref{eq3.1} satisfies
the claim, by virtue of the
$Lip(1,\frac{1}{2})$ character of $A$.

Next,
\begin{equation*}\begin{split}|\iint II
dsdt|&=|\int_{\Phi_{x+h}(E_\lambda)}
B(x,\Phi
^{-1}_{x+h}(y))f(y)J\Phi^{-1}_{x+h}
(y)dy\\ &\quad - \int_{\Phi_x(
\xi_\lambda
)}B(x,\Phi_x^{-1}(y))f(y)J\Phi^{-1}_x
(y)dy\end{split}\end{equation*} (where
$\Phi_x$ is the mapping defined in
\eqref{eq3.3})

\begin{equation*}\begin{split}&=|
\int_{\Phi_{x+h}(E_\lambda)}
\frac{[B(x,\Phi_{x+h}^{-1}(y))
-B(x,\Phi^{-1}_x(y)]}{J
\Phi_{x+h}(\Phi^{-1}_{x+h}(y)]}f(y)dy\\
&\quad+\int_{\Phi_{x+h}(E_\lambda)}
\left[\frac{1}{J\Phi_x(
\Phi^{-1}_{x+h}(y))}-
\frac{1}{J\Phi_x(
\Phi_x^{-1}(y))}\right]
B(x,\Phi_x^{-1}(y))f(y)dy\\ &\quad
+\int_{\Phi _{x+h}(E_\lambda)\backslash
\Phi_x (E_\lambda
)}\frac{B(x,\Phi_x^{-1}(y))}{J\Phi
_x(\Phi_x^{-1}(y))}f(y)dy -\int_{\Phi_x
(E_\lambda )\backslash \Phi _{x
+h}(E_\lambda )} \frac{B(x,\Phi_\lambda
^{-1}(y))}{J\Phi _x
(\Phi_x^{-1}(y))}f(y)dy|\\ &\equiv
|\int(III+IV
+V+VI)dy|,\end{split}\end{equation*}
where we have used that $J\Phi_x(s,t)$
is independent of $x$. Then

\begin{align}\label{eq3.5} |\int
IIIdy|&\leq c\|f\|_\infty
\int_{\Phi_{x+h}(E_\lambda)}|B
(x,\Phi_{x+h}^{-1}(y))-B(x,
\Phi^{-1}_x(y))|J\Phi^{-1}_{x+h} (y)dy\\
&=C\|f\|_{\infty}\int_{E_\lambda}|
B(x,s,t)-B(x,s_h,t_h)|dsdt,\notag\\
\intertext{where}
\label{eq3.6}(s_h,t_h)&\equiv \Phi_x
^{-1}(\Phi_{x+h}(s,t)),\end{align} i.e.
\begin{align*}&&x+\gamma (s_h)-\gamma
(t_h)&=x+h+\gamma (s)-\gamma (t)\\
&\Leftrightarrow &
\gamma (s_h)-\gamma (t_h)&=h+\gamma
(s)-\gamma (t)\\ &\Leftrightarrow &
(s_h-t_h,s^2_h-t^2_h)&=(h_1+s-t,
h_2+s^2-t^2).\end{align*} (by definition
of $\gamma$) i.e.,
\begin{align*}
s_h&=\frac{1}{2}\left\{h_1+s-t+
\frac{h_2+s^2-t^2}{h_1+s-t}\right\}\\
t_h&=\frac{1}{2}\left\{t-s-h_1+
\frac{h_2+s^2-t^2}{h_1+s-t}
\right\}\end{align*} Now,
\begin{equation}\begin{split}
\label{eq3.7}|B(x,s,t)-B(x,s_h,t_h)
|&\equiv |[A(x)-A(x-\gamma
(t))][A(x-\gamma (t))-A(x+\gamma
(s)-\gamma (t)]\\ &\quad-
[A(x)-A(x-\gamma (t_h))][A(x-\gamma
(t_h)-A(x+\gamma (s_h)-\gamma (t_h)]|\\
&=|[A(x)-A(x-\gamma
(t))-\{A(x)-A(x-\gamma
(t_h))\}][A(x-\gamma (t))\\ &\quad
-A(x+\gamma (s)-\gamma (t)]\\&\quad +
[A(x)-A(x-\gamma (t_h)][A(x-\gamma
(t))-A(x+\gamma (s)-\gamma (t))\\ &\quad
-\{A(x-\gamma (t_h)-A(x+\gamma
(s_h)-\gamma (t_h)\}]|\\ &\leq
C(\|\gamma (t_h)-\gamma (t)\|+\|\gamma
(s_h)-\gamma
(s)\|),\end{split}\end{equation} where
$\|\cdot \|$ denotes the parabolic metric
$\|(u,v)\|\cong |u|+|v|^{1/2}$. Notice
that
\begin{equation*}|t^2_h-t^2|\cong
|t_h-t|,\end{equation*} and similarly
for $s^2_h-s^2$, so that
\begin{align*}\|\gamma (t_h)-\gamma
(t)\|&\leq C|t_h-t|^{1/2}\\
\|\gamma (s_h)-\gamma (s)\|&\leq
C|s_h-s|^{1/2}.\end{align*} Furthermore
\begin{equation}\begin{split}\label{eq3.8}|t^\cdot
-t_h|&=\left|\frac{1}{2}
\left\{t+s+h_1-
\left(\frac{h_2+s^2-t^2}{h_1+s-t}
\right)\right\}\right|\\
&=\left|\frac{h_1}{2}+\frac{1}{2}
\left\{\frac{(h_1+(s-t))(t+s)-
(h_2+s^2-t^2)}{h_1+s-t}\right\}\right|\\
&=\left|\frac{h_1}{2}+\frac{1}{2}
\left\{\frac{h_1(t+s)-h_2}{h_1+s-t}
\right\}\right|\\ &\leq
\lambda ^{2/3},\end{split}\end{equation}
since $|h_1|\leq \lambda \leq 1$,
$|h_2|\leq \lambda ^2$, and
$s-t>15\lambda ^{1/3}$, $\frac{1}{2}\leq
s,t\leq 2$, on $E_\lambda$. A similar
estimate holds for $|s-s_h|$. Hence
$\int(III)dy\leq C\lambda ^{1/3}$ as
desired, by virtue of \eqref{eq3.5} and
\eqref{eq3.7}. Next, we observe that
\begin{equation*}\begin{split}|\int IV
dy|&\leq C\|f\|_\infty
\int_{\Phi_{x+h}(E_\lambda)}\left|
\frac{1}{J\Phi_x(\Phi_{x+h}^{-1}(y))}-
\frac{1}{J\Phi_x(\Phi_x^{-1}(y))}
\right|dy\\
&=C\|f\|_\infty\iint_{E_\lambda}
\left|\frac{1}{s-t}-\frac{1}{s_h-t_h}
\right|{(s-t)}ds dt
\end{split}\end{equation*} where
$(s_h,t_h)$ is defined as above
(see \eqref{eq3.6}), and where we have used
that $J\Phi _x(s,t)=2(s-t)$. As we have
observed (see \eqref{eq3.8}),
$|t-t_h|+|s-s_h|\leq 2\lambda
^{2/3}<<15\lambda ^{\frac{1}{3}}\leq
s-t$, so that
\begin{equation*}\begin{split}|\int IV
dy|&\leq C\|f\|_\infty
\int \int_{E_\lambda}\frac{\lambda
^{2/3}}{(s-t)^2}(s-t) ds dt\\ &\leq
C\|f\|_\infty \lambda
^{1/3}\end{split}\end{equation*} as
desired.\end{claim}

Turning last to the term $|\int Vdy|$
(the term $|\int VI dy|$ can be handled
by similar arguments, which we omit) we
see that, since $J\Phi_x(s,t)=2(s-t)$,
\begin{equation*}\begin{split}|\int
Vdy|&\leq C\|f\|_\infty
\int_{\Phi_{x+h}(E_\lambda)\backslash
\Phi_x(E_\lambda
)}\frac{1}{J\Phi_x(\Phi^{-1}_x(y))}dy\\
&\leq C\|f\|_\infty
\iint_{F_\lambda}\frac{1}{s_h-t_h}
(s-t)dsdt,\end{split}\end{equation*}
where $F_\lambda
=\Phi_{x+h}^{-1}(\Phi_{x+h}(
E_\lambda)-\Phi_x(E_\lambda))$, and where
$(s_h,t_h)$ is defined as above
(see \eqref{eq3.6}). Since $s_h-t_h\geq
c\lambda^{1/3}$, it is enough to show
that $|F_\lambda |\leq C\lambda^{2/3}$.

To this end, suppose that $(s,t)\in
F_\lambda$. Then, in particular,
$(s,t)\in E_\lambda$, and
$y\equiv
\Phi_{x+h}(s,t)\in
\Phi_{x+h}(E_\lambda)\backslash
\Phi_x(E_\lambda)$, i.e., $y\neq
\Phi_\lambda (s',t')$ for any
$(s',t')\in E_\lambda$. On the other
hand, we have observed previously that
$y=\Phi_x(s_h,t_h)$ for some
$(s_h,t_h)$ satisfying
$$|s_h-s|+|t_h-t|\leq 12\lambda
^{2/3}.$$ Since $(s_h,t_h)\notin
E_\lambda$, $(s,t)\in E_\lambda$, we
must have
$\dist ((s,t),\partial E_\lambda)\leq
C\lambda^{2/3}$. Thus $|F_\lambda |\leq
C\lambda^{2/3}$ as desired, and the
proof is done.

\end{document}